\documentclass[11pt,twoside]{article} 
\usepackage{hyperref,
amsmath,amsfonts,
graphicx,psfrag,fancybox,
qsymbols,indentfirst}
\usepackage[latin1]{inputenc}
\begin{document}
%
\newtheorem{lem}{Lemma}[section]
\newtheorem{defi}[lem]{Definition}
\newtheorem{defilem}[lem]{Definition and Lemma}
\newtheorem{theo}[lem]{Theorem}
\newtheorem{cor}[lem]{Corollary}
\newtheorem{prop}[lem]{Proposition}
\newtheorem{rem}[lem]{Remark}
\newtheorem{fig}{\hspace{2cm} Figure}
\newqsymbol{`P}{\mathbb{P}}
\newqsymbol{`Q}{\mathbb{Q}}
\newqsymbol{`E}{\mathbb{E}}
\newqsymbol{`N}{\mathbb{N}}
\newqsymbol{`R}{\mathbb{R}}

\def \imp{\Rightarrow}
\def \be{\begin{eqnarray*}}
\def \ee{\end{eqnarray*}}
\def \ben{\begin{eqnarray}}
\def \een{\end{eqnarray}}
\def\d{\displaystyle}
\def\pn{\par\noindent}
\def\mod{mod\`ele}
\def\sar{{\cal S}^{\downarrow}}
\def\mbis{M^{BIS} _n (\theta)}
\def\mbst{M^{BST}_n}
\def\proof{\noindent{\bf Proof:}\hskip10pt}
\def\ex{\noindent{\bf Examples:}\hskip10pt}
\def\QED{\hfill\vrule height 1.5ex width 1.4ex depth -.1ex \vskip20pt}
\newdimen\AAdi%
\newbox\AAbo%
\font\AAFf=cmex10
\def\AArm{\fam0 }
\def\AAk#1#2{\setbox\AAbo=\hbox{#2}\AAdi=\wd\AAbo\kern#1\AAdi{}}%
\def\AAr#1#2#3{\setbox\AAbo=\hbox{#2}\AAdi=\ht\AAbo\raise#1\AAdi\hbox{#3}
}%
\def\BBa{{\AArm A}}%
\def\BBb{{\AArm I\!B}}%
\def\BBc{{\AArm C\AAk{-1.02}{C}\AAr{.9}{I}{\AAFf\char"3F}}}%
\def\BBd{{\AArm I\!D}}%
\def\BBe{{\AArm I\!E}}%
\def\BBf{{\AArm I\!F}}%
\def\BBg{{\AArm G\AAk{-1.02}{G}\AAr{.9}{I}{\AAFf\char"3F}}}%
\def\BBh{{\AArm I\!H}}%
\def\BBi{{\AArm I\!I}}%
\def\BBj{{\AArm J}}%
\def\BBk{{\AArm I\!K}}%
\def\BBl{{\AArm I\!L}}%
\def\BBm{{\AArm I\!M}}%
\def\BBn{{\AArm I\!N}}%
\def\N{{\AArm I\!N}}%
\def\BBo{{\AArm O\AAk{-1.02}{O}\AAr{.9}{I}{\AAFf\char"3F}}}%
\def\BBp{{\AArm I\!P}}%
\def \dis{\displaystyle}
\def\P{{\AArm I\!P}}%
\def \p{{\cal P}}
\def\BBq{{\AArm Q\AAk{-1.02}{Q}\AAr{.9}{I}{\AAFf\char"3F}}}%
\def\Q{{\AArm Q\AAk{-1.02}{Q}\AAr{.9}{I}{\AAFf\char"3F}}}%
\def\BBr{{\AArm I\!R}}%
\def\BBs{{\AArm S}}%
\def\BBt{{\AArm T\AAk{-.62}{T}T}}%
\def\BBu{{\AArm U\AAk{-1}{U}\AAr{.95}{I}{\AAFf\char"3F}}}%
\def\BBv{{\AArm V}}%
\def\BBw{{\AArm W}}%
\def\BBx{{\AArm X}}%
\def\BBy{{\AArm Y}}%
\def\BBz{{\AArm Z\!\!Z}}%
\def\BBone{{\AArm 1\AAk{-.8}{I}I}}%
\def\D{{\cal D}}
\def\R{{\math R}}
\def\N{{\math N}}
\def\E{{\math E}}
\def\P{{\math P}}
\def\Z{{\math Z}}
\def\Q{{\math Q}}
\def\C{{\math C}}
\def \e{{\rm e}}
\def \f{{\cal F}}
\def \g{{\cal G}}
\def \L{{\cal L}}
\def \d{{\tt d}}
\def \k{{\tt k}}
\def\h{\frac{1}{2} \BBz^-}
\def \s{{\cal S}}
\def\u{{\upsilon}}
\newcommand{\ed}{\mbox{$ \ \stackrel{d}{=}$ }}
\def\videbox{\mathbin{\vbox{\hrule\hbox{\vrule height1ex \kern.5em\vrule height1ex}\hrule}}}
\font\calcal=cmsy10 scaled\magstep1

\def \build#1#2#3{\mathrel{\mathop{\kern 0pt#1}\limits_{#2}^{#3}}}
\def\suPP#1#2{{\displaystyle\sup _{\scriptstyle #1\atop \scriptstyle #2}}}
\def\proDD#1#2{{\displaystyle\prod _{\scriptstyle #1\atop \scriptstyle #2}}}
\def\proof{\noindent{\bf Proof:}\hskip10pt}
\def\QED{\hfill\vrule height 1.5ex width 1.4ex depth -.1ex \vskip20pt}
\setcounter{page}{1}
\thispagestyle{empty}

\markboth{~~\hrulefill~~ Chauvin and Rouault}
{Connecting tree models via martingales ~~\hrulefill~~ }

\topmargin=0mm
\vspace{3.5cm}
{\Large \bf \noindent Connecting Yule process, Bisection  and Binary Search Tree via martingales}
\\[1cm]
{\bf \noindent B. Chauvin, A. Rouault}\\ \\
{\small LAMA, B\^atiment Fermat,
Universit\'e de Versailles F-78035 Versailles. e-mail:
chauvin@math.uvsq.fr, rouault@math.uvsq.fr}\\ \\

\vspace{1cm}
{\noindent \bf Abstract.} 
We present new links between some remarkable martingales found in the study of the Binary Search Tree 
or of the bisection problem, looking at them on the probability space of a continuous time binary branching process.\\ 

\renewcommand{\thefootnote}{}
\renewcommand{\thefootnote}{\arabic{footnote}}
\setcounter{footnote}{0}

\section{Introduction}
This paper is a kind of game with martingales around the Binary Search Tree (BST) model 
(see Mahmoud \cite{Mah}). Under the random permutation model, the BST process  is an increasing (in size) sequence of binary trees $({\cal T}_n)_{n\geq 0}$ storing data, in such a way that for every integer $n$, ${\cal T}_n$ has $n+1$ leaves; the growing from time $n$ to time $n+1$ occurs by choosing uniformly a leaf and replacing it by an internal node with two leaves. In the BST we are interested 
in the profile, i.e. 
the number of leaves in each generation. 
A polynomial that encodes the profile (called the level polynomial) allows us to define
 a family of martingales $(M_n^{BST}(z) , n \geq 0)$ 
where $z$ is a positive real parameter (see Jabbour \cite{Jab1}). It is defined in Section \ref{BST}.

There are (at least) two ways of connecting the BST model to branching random walks, and to take advantage of related probabilistic methods and results.

\medskip

The first one consists in embedding the BST process into a conti- nuous time process. It is a good way to 
create independence between disjoint subtrees. The Yule process (see Athreya and Ney \cite{ANey}) is a continuous time binary branching process in which an ancestor has an exponential $1$ distributed lifetime 
and at his death, gives rise to two children with independent exponential $1$ lifetimes and so on. Define the position of an individual as its generation number, and let $Z_t$ be   
the sum of Dirac masses of positions of the population living at time $t$. 
The process $(Z_t , t \geq 0)$ is a continuous time branching random walk; call it the Yule-time process. When keeping track of the genealogical structure, call 
$\BBt_t$ the tree at time $t$ and call $(\BBt_t , t \geq 0)$ a Yule-tree process.

By embedding, the BST is a Yule-tree  process stopped at $\tau_n$, the first time when $n+1$ 
individuals exist (Pittel \cite{Pitt86}, Biggins \cite{BigMin}, Devroye \cite{Dev1}). 
A continuous time family of martingales $(M(t,z), t \geq 0)$ is attached to this model; 
it is defined in Section \ref{three}. 
In a recent paper (Chauvin \& al. \cite{ckmr}) where several models are embedded in the probability space of the  Yule-tree  process, the martingale $(M_n^{BST}(z), n \geq 0)$ appears as a  projection of the martingale $(M(t,z), t \geq 0)$.

In addition, consider the Yule tree and, on each branch, the successive birthdates of descendants of the ancestor. They are sums of independent exponentially distributed 
random variables with mean $1$,  
 so that it is natural to exchange time and space and to look at these birthdates as successive positions in a random walk. Combined with the independence between subtrees, it gives a discrete time branching random walk; call it the Yule-generation process.

\medskip

The second way consists in ``approaching'' the BST by the so-called bisection model (Devroye \cite{Dev1}, Drmota \cite{Dr2}). It is also known as  
 Kolmogorov's rock model. An object (rock)
is initially of mass one. At time $1$ it is broken into two rocks with uniform size. At time $n$, each rock (there are $2^n$) is broken independently from the other ones into two rocks with uniform size.
The mass of each rock results from the product of independent uniform random variables.  
Taking logarithms gives a discrete time branching random walk; call it the Bisection process.  
In \cite{ckmr}, 
it was observed in the Yule-tree environment.

\medskip

For these reasons, it is worth  considering all these models on the same probability space. As explained in Section 2, there are three branching random walks: Yule-time, Yule-generation and Bisection, each one with its family of additive martingales. 

In Section 3, thanks to a convenient adjustment of parameters of these four families (the three previous ones and the BST), we establish 
strong links between these martingales and their limits as $n$ (or $t$) tends to infinity. 
In Theorems \ref{bis=bst} and \ref{MGEN=YULE}, we claim that
there are actually 
only two different limits (a.s.) in the domain of $L^1$ convergence. On the boundary of this domain, 
we identify in Theorem \ref{derivMGEN=YULE} limits of martingales obtained by taking 
derivatives  with respect to the parameter.   

To prove these identifications, we need uniqueness arguments which are explained in Section 4. 
We write the (stochastic) equations satisfied by 
the limits of martingales. The solutions of these equations have distributions which
are fixed points of so-called smoothing transformations, as defined in Holley and Liggett \cite{HolLig} 
or in  Durrett and Liggett
\cite{DurLig}. Owing to known results on uniqueness of their Laplace transforms
(Liu \cite{Liu98}, \cite{Liu00}, Kyprianou \cite{Kyp98}, Biggins and Kyprianou \cite{BK1}) 
we  get equalities in law (Proposition 4.1).

Section 5 is devoted to proofs of theorems of Section 3. 
In particular, we show that equalities in law are (a.s.) equalities between random variables. 

In Section 6, we explain some relations between the above 
 functional equations satisfied by  Laplace transforms and equations 
 studied by Drmota in recent papers on the height of the BST (\cite{Dr0,Dr1}). This allows to get 
 explicit limiting distributions in the case $z= 1/4$.

\medskip

Let us now fix some notation. 
In the whole paper we are concerned with binary trees whose nodes (also called individuals) are labelled by the elements of
\[\mathbb{U} :=\{\emptyset\}\cup\bigcup_{n\geq 1}\{0,1\}^n\,,
\]
the set of finite  words on the alphabet $\{0,1\}$ (with $\emptyset$ as an empty word). For $u$ and $v$ in $\mathbb{U}$, denote by $uv$ the concatenation of the word $u$ with the word $v$ (by convention we set, for any $u\in\mathbb{U}$, $\emptyset u=u$). 
If $v \not= \emptyset$, we say that $uv$ is a descendant of $u$ and $u$ is an ancestor of $uv$, in particular $v$ is the father of $v0$ and $v1$. We note $u\succ v$ to say that $v$ is an ancestor of $u$.
A \sl complete binary tree \rm $T$ is a finite subset  of $\mathbb{U}$ such that
\[\left\{
\begin {array}{l}
\emptyset\in T\\
\textrm{ if }uv\in T \textrm{ then } u\in T \,,\\
u1 \in T \Leftrightarrow u0\in T\,.
\end {array}
\right.\]
The elements of $T$ are called \it nodes \rm, and $\emptyset$ is called the \it root \rm ; \rm $|u|$, the number of letters in $u$, is the \it depth \rm of $u$ (with $|\emptyset|=0$).
Write {\bf BinTree} for the set  of complete binary trees.

A tree $T \in$ {\bf BinTree} can be described by giving the set $\partial T$ of its \it  leaves\rm, that is, the nodes that are  in $T$ but with no descendants in $T$.
The nodes of $T\backslash\partial T$ are called \it internal \rm nodes.

\vspace{0.3cm}
\section
{The four martingales}
\label{martingales}
\vspace{0.3cm}
\subsection{Branching random walks}
\label{brw}

A discrete time supercritical branching random walk (in $\BBr$) is recursively defined as follows: the initial ancestor is at the origin 
 and the positions of his children form a point process $Z$. The distribution of this point process $Z$ is a probability 
 on $M$, the set of locally finite sums of Dirac measures. Each child of the ancestor reproduces in the same way and each individual also does: the positions of each sibling relative to its parent are an independent copy of $Z$. Let $Z_n$ be the point process in $\BBr$ formed by the $n$-th generation. The intensity of $Z$ is the Radon measure $\mu$ defined for every nonnegative bounded function $f$ by 
$$\int_{\BBr} f(x) \mu(dx) = `E\big(\int_{\BBr} f(x) Z(dx)\big) \,.$$
and the intensity of $Z_n$ is $\mu^{*n}$.
We assume $1 < \mu(\BBr) \leq +\infty$ (supercriticality).  

We define for $\theta \in \BBr$
\be
\Lambda (\theta) = \log `E\int_\BBr e^{\theta x} Z(dx)\,.
\ee
The (positive) martingale associated with this process is   
$${\bf M}_n (\theta)= \int_\BBr e^{\theta x - n \Lambda (\theta)} Z_n (dx)\,.$$
Let ${\bf M}_\infty(\theta)$ be the a.s. limit.
Under a ``$k \log k$'' type condition, we have (Biggins' convergence theorem, for instance in \cite{JDB77,JDB92})

- if $\theta \Lambda'(\theta) - \Lambda (\theta) < 0$,  the convergence is also in $L^1$ and $`E{\bf M}_\infty(\theta) = 1$,

- if $\theta \Lambda'(\theta) - \Lambda (\theta) \geq  0$, then ${\bf M}_\infty(\theta) = 0$ a.s.

By analogy with the Galton-Watson process, we call ``supercritical'' the values of $\theta$ in the first region, 
and ``critical'' (resp. ``subcri- tical'') if they correspond to equality (resp. strict inequality) in the second region.
 
In a continuous time branching random walk (in $\BBr$), the starting point is the same as above. Each  individual has an independent exponential lifetime (of
parameter $\beta$), does not move during its life, and at its death is replaced by children according to a copy of a point process $Z$ exactly as in the discrete-time scheme. 
The role of $\Lambda$ is now played by 
$$L(\theta) = \beta \big( `E\int_\BBr e^{\theta x} Z(dx) - 1\big)\,.$$
If $Z_t$ denotes the random measure of positions of individuals alive at time $t$, the (positive) martingale associated to this process is   
$${\bf M} ( t, \theta)= \int_\BBr e^{\theta x - t L (\theta)} Z_t (dx)\,.$$
Its behavior as $t \rightarrow \infty$ is similar to the above, with $L$ instead of $\Lambda$ (Biggins \cite{JDB92}, Uchiyama \cite{Uchi}). We denote by ${\bf M} ( \infty, \theta)$ its limit.

At critical values of $\theta$, the limit martingale vanishes. It is then classical, since Neveu \cite{Nev1}, 
to study the family of derivatives
$
\frac{\partial}{\partial \theta} {\bf M}_n (\theta)
$  (or  $\frac{\partial}{\partial \theta} {\bf M} (t, \theta)$). 
These are martingales of expectation $0$, 
which converge a.s. to 
a random variable of constant sign, of infinite expectation, under appropriate conditions (Kyprianou \cite{Kyp98}, Liu \cite{Liu00}, Biggins-Kyprianou \cite{BK1} and Bertoin-Rouault \cite{BertRou1}). The  details are given below.

\vspace{0.3cm}
\subsection{Examples}\label{three}

\noindent{\bf{Bisection martingale}}

Passing to logarithms in the bisection problem, we get 
 a discrete time branching random walk 
whose
 reproduction measure is $ Z^{BIS}= \delta_{-\log U} + \delta_{-\log (1-U)}$, where $U\sim {\cal U}([0,1])$\footnote{${\cal U}([0,1])$ is the uniform distribution on $[0,1]$,
${\cal E}(\lambda)$ is the exponential distribution of parameter 
$\lambda$,  $\sim$ means ``is distributed as'' and $\buildrel{law}\over{=}$ means equality in distribution.}. 
We have $$ \Lambda (\theta) = \log 
`E [U^{-\theta} + (1-U)^{-\theta}] 
= \log \frac{2}{1- \theta}\,.$$
Let us make a change of parameter, setting
$z = \frac{1 - \theta}{2}$, 
 so that 
$\Lambda(\theta) = -\log z$. The corresponding martingale is
\ben
\label{defbis}
M^{BIS}_n (z) := \sum_{|u| = n} \  e^{(1 -2z)X_u} z^n 
\,,
\een
where $X_u$ denotes the position of individual $u$. 
It is easy to see that the range of $L^1$ convergence is $z \in (z_c ^- , z_c ^+)$, where $z_c^- < z_c ^+$ are the two (positive) solutions of
\be
 2z \log {z} - 2z + 1 = 0\,,
\ee
(i.e. $c'$ and $c$ in Drmota's notation, see \cite{Dr1})
\be
z_c ^- = \frac{c'}{2} = 0.186\ldots , \ \ z_c ^+ = \frac{c}{2} = 2.155\ldots 
\ee 
For $z = z_c ^-$ (resp. $z_c ^+$), applying Theorem 2.5 of Liu \cite{Liu00}, we see that 
the derivative $(M_n^{BIS})' (z)$ converges a.s. to a limit denoted by $M'^{BIS}_\infty (z)$ which is positive (resp. negative) and has infinite expectation.

\medskip
 
\noindent{\bf{Yule-time martingale}} 

The Yule-time process is a continuous time branching random walk, its reproduction measure is $Z=2\delta_{-\log 2} $, and the parameter $\beta$ of the exponential lifetimes is equal to $1$. We have $L( \theta)=2^{1-\theta}-1$. The position of individual $u$ at time $t$ is
 $X_u(t) = -|u| \log 2$.
Introducing the parameter $z=2^{-\theta}$ we have $L(\theta) = 2z -1$ and the corresponding martingale becomes
\be
M(t,z) = \sum_{u \in 
{\cal Z}_t} z^{|u|} e^{t(1 -2z)}
\ee
where ${\cal Z}_t$ denotes the set of individuals alive at time $t$, of cardinality $N_t$. This can be considered as a generalization of the classical Yule martingale $(e^{-t} N_t , t\geq 0)$ which is known to converge a.s. to a random variable $\xi \sim {\cal E}(1)$. The behavior of this family follows the same rule as above.
Moreover, it was proved in Bertoin-Rouault \cite{BertRou1} (see also \cite{ckmr}) that for $z =z_c^\pm$, the derivative $M'(t,z)$ converges a.s. to  a limit denoted by $M'(\infty, z)$, of constant sign and infinite expectation.

\medskip

\noindent{\bf{Yule-generation martingale}} 

The Yule-generation process is a discrete time branching random walk and its reproduction measure is $Z = 2 \delta_{\epsilon}$ where $\epsilon \buildrel{law}\over{=}
 {\cal E}(1)$ and the factor $2$  
means that the two brothers  
 appear at the same time. Since $\epsilon \buildrel{law}\over{=}
 -\log U $, the intensity $\mu$ is the same as in the bisection case, and then 
we have  $\Lambda (\theta) = \log \frac{2}{1-\theta}$ again. 
With the same change of parameter, we have a martingale
\ben
\label{defmgen}
M^{GEN}_n (z) =  \sum_{|u| = n} \  e^{(1 -2z)X_u} z^n 
\,,
\een which has the same form as $M^{BIS}_n (z)$, and has the same range of $L^1$ convergence. 
 However the martingales $M^{BIS}_n (z)$ and $M^{GEN}_n (z)$ do not have the same distribution 
since the dependence between the positions $X_u , |u| = n$, is different in the two models 
(although the structure of random variables along a given branch is the same).
Again, from Liu \cite{Liu00}, we have for $z = z_c^\pm$ convergence a.s. of  
$(M_n ^{GEN})'(z)$ to a limit $M'^{GEN}_\infty(z)$, of constant sign and infinite expectation.

\vspace{0.3cm}
\subsection{The BST martingale}\label{BST}

A binary search tree (BST) process (for a detailed description, see Mahmoud \cite{Mah}) is a sequence $ ({\cal T}_n , n \geq 0)$ of complete binary trees, where ${\cal T}_n$ has $n$ internal nodes, which grows by successive insertions of data, 
under the so-called random permutation model.
Let us describe the dynamics of the sequence of trees. Tree ${\cal T}_{1}$ is reduced to the root and has two leaves. Tree ${\cal T}_{n+1}$ is obtained from ${\cal T}_n$ 
by replacing one of its $n+1$ leaves by an internal node and thus creating two new leaves. The insertion is done uniformly on the leaves, which means with probability $1/(n+1)$.

To study the shape of these trees, it is usual to 
 define the profile of tree ${\cal T}_n$ by the collection of
\be
U_k(n):= \# \{ u \in \partial {\cal T}_n , |u| = k \} \ \ , \ \ k\geq 1\,,
\ee
counting the number of leaves of ${\cal T}_n$ at each level. The profile is coded by the level polynomial $\sum_k U_k (n) z^k$, for $ z \in [0, \infty)$ and can be studied by martingale methods (Jabbour \cite{Jab1}, Chauvin \& al. 
\cite{Jab2,ckmr}).
  Because of the dynamics of the tree process,  this polynomial, renormalized by its expectation, is a ${\cal F}_{(n)}$-martingale, where ${\cal F}_{(n)}$ is the $\sigma$-field generated by all the events $\{ u \in {\cal T}_j\}_{j \leq n , u \in \mathbb{U}}$. More precisely, we define the BST martingale
\be
\mbst (z) = \sum_{u \in \partial {\cal T}_n} \frac{z^{|u|}}{C_n (z)} =  \sum_k  \frac{U_k (n)}{C_n (z)}z^k
\ee
where $C_0(z)=1$ and
\be
C_n (z) = \prod_{k=0}^{n-1} \frac{k + 2z}{k+1} = (-1)^n 
\begin{pmatrix} - 2z\\n\end{pmatrix} , \hskip 0.5cm \  n\geq 1\,. 
\ee

In the supercritical range $z \in (z_c ^- , z_c ^+)$, this martingale converges in $L^1$ to a nondegenerate limit $M^{BST}_\infty (z)$ and converges a.s. to $0$ elsewhere, in particular for the critical values $z_c ^-$ and $z_c ^+$. For these critical values, the derivative $(M_n^{BST})' (z)$ is a martingale of expectation $0$, which converges a.s. to  a random variable $M'^{BST}_\infty (z)$ of constant sign, of infinite expectation.

\section{Connections between these martingales}\label{connections}

We now set all these martingales on the same probability space and we define below the continuous time tree valued Yule process.  Roughly speaking, 

- in the BST we keep track of profile, 

- in the Bisection, we keep track of the balance between right subtrees and left subtrees,

- in the Yule-generation, we keep track of time of appearance of the different nodes.  
 
This set-up provides nice connections which are precised in Theo- rems \ref{bis=bst}, \ref{MGEN=YULE} and \ref{derivMGEN=YULE}. 
\medskip

Let $({\u}_t)_{t\geq 0}$ be a Poisson point process taking values in $\mathbb{U}$ with 
intensity measure  $\nu_{\mathbb{U}}$, the counting measure on $\mathbb{U}$.
 Let $(\BBt_t)_{t \geq 0}$  be
a  {\bf BinTree} valued process starting from
 $\BBt_0 = \{\emptyset \}$ and  jumping only when $({\u}_t)_{t\geq 0}$  jumps. 
Let $t$ be a jump time for $\u_\cdot$; $\BBt_t$ is obtained from $\BBt_{t-}$ in the following way:

if $\u_t \notin \partial \BBt_{t-}$ keep $\BBt_t = \BBt_{t-}$ and if $\u_t \in \partial \BBt_{t-}$ 
take $\BBt_t = \BBt_{t-}\cup \{\u_t0 , \u_t1\}$.

\noindent The counting process $(N_t)_{t \geq 0}$ defined by
\be
N_t := \# \partial\BBt_t
\ee
is the classical Yule (or binary fission) process (Athreya-Ney \cite{ANey}).
In the following, we refer to the continuous-time tree process $(\BBt_t)_{t \geq 0}$ as the Yule tree process.

We note $0=\tau_0<\tau_1<\tau_2<\ldots$  the successive jump times (of $\BBt.$),
\be
\tau_n = \inf \{t : N_t = n+1\}\,.
\ee

\begin{figure}[htbp]
\psfrag{a}{0}
\psfrag{b}{$\tau_1$}
\psfrag{c}{$\tau_2$}
\psfrag{d}{$\tau_3$}
\psfrag{e}{$\tau_4$}
\psfrag{f}{$\tau_5$}
\psfrag{g}{$t$}
\centerline{\includegraphics[height=5cm]{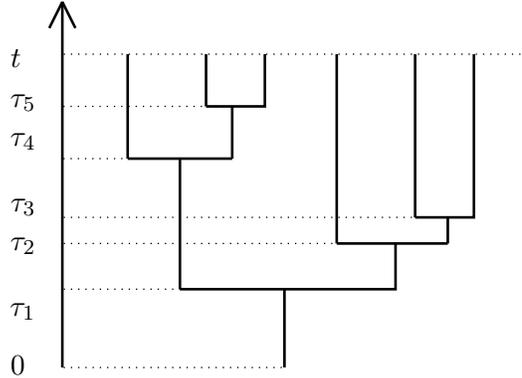}}
\caption{\label{frag2} Continuous time binary branching process. }
\end{figure}

We define recursively the time of appearance (or time of saturation) of nodes  by
$$S^\emptyset = 0 , \ \ S^{u_1 .. u_n} = \inf \{ s > S^{u_1 .. u_{n-1}} : \u_s = u_1 .. u_{n}\}$$
(the definition of $\u_t$ is given above). Actually $S^{u_1 .. u_n}$ 
is the sum of $n$ i.i.d. ${\cal E}(1)$ random variables. This yields
$\BBt_t = \{ u : S^u \leq t\}$. The natural filtration is $({\cal F}_t , t\geq 0)$ where ${\cal F}_t$ is 
generated by all the random variables $\u_s, s \leq t$.
Another useful one is 
$({\cal F}^n , n \geq 1)$ where ${\cal F}^n$ is generated by the variables $S^v$ for all $|v| \leq n$. 
Finally we will use $({\cal F}_{(n)} , n \geq 1)$, 
where ${\cal F}_{(n)}$ is generated by ${\BBt}_{\tau_1}, \ldots ,{\BBt}_{\tau_n}$.

\medskip
 
This Yule-time process can also be seen as a fragmentation process.
We may encode dyadic open subintervals of $[0,1]$ with elements of $\mathbb{U}$: we set $I_\emptyset = (0,1)$ and for $u = u_1 u_2 \ldots u_k \in \mathbb{U}$,
$$I_u = \Big(\sum_{j=1}^k u_j 2^{-j}, 2^{-k}+\sum_{j=1}^k u_j 2^{-j} \Big).$$
With this coding, the  evolution corresponding to the Yule-time process is a very simple example of 
fragmentation process. This idea goes back to Aldous and Shields (\cite{AlS} Section 7f and 7g).
In other words, for $t\geq 0$, $F(t)$ is a finite family of intervals. At time 0, we have $F(0)=(0,1)$. Identically independent exponential ${\cal E}(1)$ random variables
  are associated with intervals of $F(t)$.
Each interval in $F(t)$ splits into two parts (with same size)  independently  of each other after an exponential time ${\cal E}(1)$.\par
Hence, one has $F(0)=(0,1)$, $F(\tau_1)=((0,1/2),(1/2,1))$ where $\tau_1\sim {\cal E}(1)$, etc... 
One can interpret the two fragments $I_{u0}$ and $I_{u1}$ issued from $I_u$ as the two children of $I_u$, one being the left (resp. right) fragment $I_{u0}$ (resp. $I_{u1}$), obtaining thus a binary tree structure. 
An interval with length $2^{-k}$ corresponds to a leaf at depth $k$ in the corresponding tree structure.
\begin{figure}[htbp]
\psfrag{a}{0}
\psfrag{b}{$\tau_1$}
\psfrag{c}{$\tau_2$}
\psfrag{d}{$\tau_3$}
\psfrag{e}{$\tau_4$}
\psfrag{f}{$\tau_5$}
\psfrag{g}{$t$}
\psfrag{0}{0}
\psfrag{1}{1}
\centerline{\includegraphics[height=6cm]{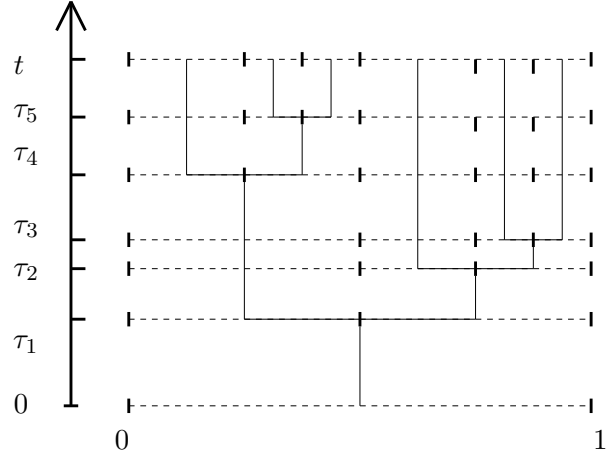}}
\caption{\label{frag} Fragmentation and its tree representation.}
\end{figure}
\medskip

\vspace{0.3cm}
\subsection{Connection Yule-time $\rightarrow $ BST}

In  \cite{ckmr}, it is proved that 
\be
\bigl(\BBt_{\tau_n} , n \geq 1\bigr) \buildrel{law}\over{=} \bigl( {\cal T}_n , n \geq 1\bigr)\,.
\ee

We can now consider the BST process as the Yule process observed at the (random) splitting times $\tau_n$. It turns out that 
$M^{BST}_n (z)$ is the projection of $M(\tau_n , z)$ on ${\cal F}_{(n)}$.
It yields nice limit martingale connections    
\ben
\label{lmc1}
M(\infty, z) &=& \frac{\xi^{2z-1}}{\Gamma(2z)}M^{BST}_\infty (z)\ \ \ \hbox{for} \ \ z \in (z_c^- , z_c^+)\,, \ \ a.s.\\
M'(\infty,z)  &=& \frac{\xi^{2z-1}}{\Gamma(2z)} M'^{BST}_\infty (z)\ \ \ \hbox{for} \ \ z = z_c^\pm\,, \ \ a.s.\,.
 \een
where $\xi \sim {\cal E}(1)$.

\vspace{0.3cm}
\subsection{Connection BST  $\rightarrow $ Bisection}

For $v \in \mathbb{U}$, let $n^{(v)}_t$ be the number of individuals alive at time $t$ in the subtree beginning at node $v$. 
In \cite{ckmr} Section 2.5, it is shown that the random variables 
\ben
\label{defU}
U^{(v)} :=  \lim_t \frac{n_t ^{(v)}}{n_t^{(father \ of \ v)}}
\een
are ${\cal U}([0,1])$ distributed, 
independent along a branch and that $U^{(v0)} + U^{(v1)} = 1$ for every $v$.

It is then clear that we may construct a Bisection branching random walk with these variables. Let us call $M_n ^{BIS} (z)$ the associated martingale:
\begin{equation}
\label{redefbis}
M_n ^{BIS} (z) := \sum_{|u| = n} \big(\prod_{u\succ v} U^{(v)} \big)^{2z - 1} z^n\,.
\end{equation}
Notice that with the notation of Section 2.2, $X_u = \sum_{v \prec u} -\log U^{(v)}$.

\begin{theo}\label{bis=bst}
For $z \in (z_c^- , z_c^+)$, and $n \geq 1$, a.s. 
\begin{equation}
\label{a}
 `E[M^{BST}_\infty (z) \ | \ {\cal F}^{BIS}_n]
= M_n ^{BIS} (z)\,,
\end{equation}
and
\begin{equation}
\label{b}
M_\infty ^{BIS}(z) = M^{BST}_\infty (z)\,.
\end{equation}
\end{theo}

\begin{theo}
\label{cornew}
For  $z = z_c^\pm$, we have
$$M_\infty^{'BST}(z) \buildrel{law}\over{=} M_\infty ^{'BIS} (z)\,.$$
\end{theo}
\medskip

It is of course tempting to conjecture equality of the above random variables. 

\subsection{Connection Yule-time  $\rightarrow $ Yule-generation}
\label{simplette}

Let ${\cal L}_n := \{u : |u| =n \}$  be the set of the nodes in the $n$-th generation 
of the Yule tree process. For $u \in \mathbb{U}$, by definition of the Yule-generation process, 
position $X_u$ can be also seen as the time $S^{u}$ of appearance of node $u$, 
i.e.  the sum of the i.i.d. ${\cal E}(1)$ lifetimes  along the branch from the root to $u$. So
\ben
\label{defmgenagain}
M^{GEN}_n (z) :=  \sum_{u \in {\cal L}_n} \  e^{(1 -2z)S^{u}} z^n\,.
\een
The following theorem is analogous to Theorem \ref{bis=bst}, it is valid in the supercritical case. Theorem \ref{derivMGEN=YULE} concerns the critical case.

\begin{theo}
\label{MGEN=YULE}

For $z \in (z_c ^-
, z_c ^+)$, a.s.
\ben\label{yuleconditioned}
`E[M(\infty, z) | {\cal F}^n] = M^{GEN}_n (z)\,,
\een
and consequently
\ben\label{limitmgen}
M^{GEN}_\infty (z) = M(\infty, z)\,.
\een
\end{theo}

\begin{theo}\label{derivMGEN=YULE}
For $z = z_c^\pm$, a.s. 
\ben\label{limitderivmgen}
M'^{GEN}_\infty (z) = M'(\infty, z)\,.
\een
\end{theo}

\begin{rem}
\end{rem} A set $L$ of nodes is usually said to have the line property if no node of $L$ is an ancestor of another node of $L$. In other words, the subtrees starting from nodes of $L$ are disjoint trees.
Theorems \ref{MGEN=YULE} and \ref{derivMGEN=YULE} could 
appear as a consequence, for the particular lines ${\cal Z}_t$ and ${\cal L}_n$,  of a more  general theorem which would be: additive martingales associated to a sequence of ``lines'' tending to infinity have the same limit, independently of the choice of this sequence. This theorem holds without any serious difficulty as soon as the notion of ``line'' is precisely defined, which is necessary since several notions  exist\footnote{let us also mention the close notion of ``cutset'' in Peres \cite{Peres97}}; 
``optional lines'' in Jagers \cite{Jagers89} 
require
 measurability of the stopping  rule with respect to the {\it process} until the line. More restrictively, ``stopping lines'' in Chauvin \cite{Chau91} and Kyprianou \cite{Kyp99,Kyp00} or frosts (in the fragmentation frame, cf Bertoin \cite{Bertoin02}) 
require
 measurability of the stopping rule with respect to the {\it branch} from the root to some node of the line. The above mentioned general theorem holds for stopping lines and {\it not} for optional lines in the Jagers' sense. Of course ${\cal Z}_t$ and ${\cal L}_n$ are stopping lines, but the stopping time 
$\tau_n$ (the first time when $n$ intervals exist in the fragmentation) defines an optional but {\it not}  stopping line.

In other words, in view of equalities $M_{\infty}^{BST}(z) =M_{\infty}^{BIS}(z) $ and $M(\infty,z)=M_{\infty}^{GEN}(z) $ and in view of connections (\ref{lmc1}) it is clear that (for $z \not= 1/2$) $M(\infty,z)$ is different from $M_{\infty}^{BST}(z)$. This is consistent with the fact that $\tau_n$ does not define a stopping line.

\vspace{0.3cm}
\section{Smoothing transformations and limit distributions}\label{distributions}

Before proving almost sure equalities announced in the theorems of  Section 3, we first look at equalities in distribution. 
The random limits mentioned in Section 2 satisfy ``duplication'' relations, which come from the binary branching structure of the underlying processes. These relations may be viewed as equalities between random variables or as functional 
equations on their Laplace transforms. The corresponding distributions are  
fixed points of  so-called smoothing transformations (Holley and Liggett \cite{HolLig}, Durrett and Liggett \cite{DurLig}).

\vspace{0.3cm}
\subsection{Duplication relations}

1) Let us first consider $M(\infty, z)$. 
After conditioning on  the first splitting time $\tau_1$  of the Yule-tree process (recall that $\tau_1 \sim {\cal E}(1)$), 
we get (\cite{ckmr} Section 3.1), for $z \in (z_c^- , z_c^+)$ 
\ben
\label{split}
M(\infty, z) = ze^{(1-2z)\tau_1} \left(M_0 (\infty, z) + M_1 (\infty, z)\right)
\ \ \ a.s.,
\een
where the random variables $M_0 (\infty, z)$ and $M_1 (\infty, z)$ are independent, distributed as $M (\infty, z)$ 
and independent of $\tau_1$. Moreover $`P(M(\infty, z)>0)=1$. Iterating (\ref{split}) we get
\begin{equation}
\label{iter}
M(\infty, z) =z^n \sum_{|u| = n} e^{(1-2z)S^u} M_u (\infty, z) \ \ a.s.
\end{equation} 
For $z= z_c^\pm$, $M(\infty , z) = 0$ a.s. and  the relation satisfied by 
$M'(\infty, z)$ 
is the same as (\ref{split}) (mutatis mutandis).
Moreover $`P(M'(\infty, z_c ^-) > 0) =  `P((M'(\infty, z_c ^+) < 0)=1$.
\medskip

2) By definition of the Yule-generation process (\ref{defmgenagain}), we have, 
conditioning  upon the first generation, 
\ben
\label{mgen}
M_\infty^{GEN}( z) = ze^{(1-2z)\tau_1} \left(M_{0, \infty}^{GEN}( z) + M_{1, \infty}^{GEN}( z)\right)
\ \ \ a.s.,
\een
which is exactly the same equation as (\ref{split}). The same result holds for derivatives at $z = z_c^\pm$.
\medskip

3) Let us see now what happens for the BST martingale limit. 
By embedding, it is shown in \cite{ckmr} Section 3.1 that
\ben
\label{masterm}
M^{BST}_{\infty} (z) = z \left(U^{2z-1} M^{BST}_{\infty, (0)} (z) 
+ (1-U)^{2z-1}M^{BST}_{\infty, (1)} (z)\right) \ \ a.s.,
\een
where $U  \sim {\cal U}([0,1])$ is nothing but $U^{(0)}$ as defined in (\ref{defU}), where $M^{BST}_{\infty ,(0)} (z), M^{BST}_{\infty, (1)} (z)$ are independent (and independent of $U$) 
and distributed as $M^{BST}_{\infty} (z)$. For $z = z_c^\pm$ the relation is the same with $M'^{BST}_{\infty} (z)$ instead of $M^{BST}_{\infty} (z)$.

Iterating (\ref{masterm}), we get
\begin{equation}
\label{iterating}
M^{BST}_{\infty} (z) =  z^n \sum_{|u| = n} \big(\prod_{v\prec u} U^{(v)} \big)^{2z - 1}  
M^{BST}_{\infty, u} (z) \ \ a.s.
\end{equation}

4) By definition of the bisection (\ref{redefbis})
$$M_n^{BIS}(z) = \sum_{|u| = n} \big(\prod_{v\prec u} U^{(v)} \big)^{2z - 1} z^n \ \ a.s.$$
We have, conditioning upon the first generation of this process, 
\begin{equation}
\label{masterbis}
M^{BIS}_{\infty} (z) = z \left(U^{2z-1} M^{BIS}_{\infty, (0)} (z) 
+ (1-U)^{2z-1}M^{BIS}_{\infty, (1)} (z)\right) \ \ a.s.
\end{equation}
which is exactly the same equation as (\ref{masterm}).

At this stage, we see two packages ($M(\infty, z), M^{GEN}_\infty (z)$) 
and ($M^{BST}_\infty (z), M^{BIS}_\infty (z)$), which are consistent with the results of 
Theorems \ref{bis=bst} and \ref{MGEN=YULE}, but which do not yet give 
equality in law, 
owing to  lack of uniqueness. 

\vspace{0.3cm}
\subsection{Functional equations and equalities in law}\label{functional}
 The relations between random variables in the previous subsection 
 imply distributional equations which can be viewed as 
 functional equations on their Laplace transforms.

Set 
\be
j(z, x) &:=& `E \exp -x M(\infty, z) \,, \hspace{6mm}\hbox{for} \ z \in (z_c^- , z_c^+)\,,\\
j(z, x) &:=& `E \exp -x |M'(\infty, z)| \,, \ \ \hbox{for} \ z=z_c^\pm\,,\\
j^{GEN}(z, x) &:=& `E \exp -x M_\infty^{GEN}(z) \,, \hspace{4mm} \hbox{for} \ z \in (z_c^- , z_c^+)\,,\\
j^{GEN}(z, x) &:=& `E \exp -x |M_\infty^{'GEN}(z)| \,, \hspace{1mm} \hbox{for} \ z=z_c^\pm\,,\\
j^{BST}(z, x)  &:=& `E \exp -x M_\infty^{BST}( z) \,, \hspace{5mm}  \hbox{for} \ z \in (z_c^- , z_c^+)\,,\\
j^{BIS}(z, x) &:=& `E \exp -x M_\infty^{BIS}(z) \,, \hspace{5mm} \hbox{for} \ z \in (z_c^- , z_c^+)\,.
\ee

Let us summarize some results on fixed points of smoothing transformations that are needed for our study. There is a broad literature on this topic. 
One of the more recent contributions is in Biggins and Kyprianou \cite{BK1}.
We choose to give these results under the assumptions of Liu \cite{Liu98}, \cite{Liu00} (see also Kyprianou \cite{Kyp98}), which are fulfilled in our examples. 

Let us consider a branching random walk as in Section 2.1 with $Z = \sum_{i=1}^N \delta_{x_i}$
and $P(N= 0) = 0$. We defined the martingale
$${\bf M}_n (\theta) = \int_\BBr e^{\theta x - n \Lambda (\theta)} Z_n(dx)\,,$$
and the derivative $${\bf M}'_n (\theta)
= \int_\BBr ( x - n \Lambda' (\theta)) e^{\theta x - n \Lambda (\theta)} Z_n(dx)\,.$$
Let us assume 
$`E N^{1+\delta} < \infty$ 
  and $`E [{\bf M}_1 (\theta)]^{1+\delta} < \infty$,  
for some $\delta > 0$.

In the supercritical range, i.e. if $\theta\Lambda'(\theta) -\Lambda(\theta) < 0$, 
the Laplace transform 
${\bf J}(s)= `E e^{-s {\bf M}_\infty (\theta)}$ satisfies 
\ben
\label{smooth0}
{\bf J}(s) = `E \prod_{i=1}^N {\bf J}(s e^{\theta x_i - \Lambda (\theta)} )
\een
(branching property at the first splitting time)
and
\ben
\label{superliu0}
\lim_{x \rightarrow 0+} \frac{1 - {\bf J}(x)} {x} = 1
\een
($L^1$ convergence of the martingale).

Moreover for every $K > 0$ there is only one solution of 
 (\ref{smooth0}) in the class of Laplace transforms of 
 nonnegative (non degenerate)  random variables satisfying 
\ben
\label{liuk0}
\lim_{x \rightarrow 0+} \frac{1 - {\bf J}(x)} {x} = K\,.
\een

For critical values, i.e. for $\theta$ such that $\theta\Lambda'(\theta) - \Lambda (\theta) = 0$, 
the a.s. limit ${\bf M}'_\infty (\theta)$  is positive if $\theta < 0$ and negative if $\theta > 0$. 
Its Laplace transform ${\bf J}(s)= `E\, e^{- s |{\bf M}'_\infty (\theta)|}$ 
satisfies 
\ben
\label{smooth}
{\bf J}(s) = `E \prod_{i=1}^N {\bf J}(s e^{\theta x_i - \Lambda (\theta)} )\,,
\een
and
\ben
\label{superliu}
\lim_{x \rightarrow 0+} \frac{1 - {\bf J}(x)} {x |\log x|} = |\theta|^{-1}\,,
\een
(see Theorem 2.5 a) of Liu \cite{Liu98} with a slight change of notation).
Moreover for every $K >0$  there is only one solution of 
 (\ref{smooth}) in the class of Laplace transforms of nonnegative (non degenerate) 
 random variables satisfying 
\ben
\label{liuk}
\lim_{x \rightarrow 0+} \frac{1 - {\bf J}(x)} {x |\log x|} = K\,.
\een
In our setting, this yields the following identities.
\begin{prop} 

\begin{itemize}
\item[a)] For $z \in (z_c^- , z_c^+)$, we have $j(z, \cdot) = j^{GEN}(z, \cdot)$ or equivalently
\ben
\label{lemme1}
M(\infty, z) \buildrel{law}\over{=} M_\infty ^{GEN} (z)\,.
\een
\item[b)] For $z \in (z_c^- , z_c^+)$, we have $j^{BST}(z, \cdot) = j^{BIS} (z, \cdot)$ or equivalently,
\ben
\label{lemme2}
M_\infty^{BST}(z) \buildrel{law}\over{=} M_\infty ^{BIS} (z)\,.
\een
\item[c)] For critical $z$ we have
\ben
\label{contraintedc}
\lim_{x\downarrow 0} \frac{1 - j^{GEN}(z_c^+ ,x)}{x |\log x|} = \lim_{x\downarrow 0} 
\frac{1 - j^{BIS}(z_c^+ ,x)}{x |\log x|} = 
\frac{2}{2z_c^+ - 1}\,,\\
\label{contraintedc+}
\lim_{x\downarrow 0} \frac{1 - j^{GEN}(z_c^- ,x)}{x |\log x|} = \lim_{x\downarrow 0} \frac{1 - j^{BIS}(z_c^- ,x)}{x |\log x|}  = \frac{2}{1 - 2z_c^-}
\,.
\een
\end{itemize}
\end{prop}

\proof Recall (see end of Section 2.2) that the branching random walks BIS and GEN are 
different but have the same $\Lambda(\theta)$, giving then 
two versions of (\ref{smooth}). They
share the same critical points (recall  the correspondence $\theta = 1-2z$).

a) From (\ref{split}) and (\ref{mgen}) it is  clear that for $z \in [z_c^- , z_c^+ ]$, the functions $j(z, \cdot)$ and $j^{GEN}(z, \cdot)$ are  non-constant solutions of
\begin{eqnarray}
\label{eqj}
{\bf J}(x) &=& \int_0 ^1 {\bf J}(zxu^{2z-1})^2 du\,,
\\
\label{ci}
{\bf J}(0) &=& 1\,.
\end{eqnarray}
Moreover, since $`EM(\infty ,z) = `EM_\infty ^{GEN} (z) = 1$ for $z$ supercritical, it turns out that 
they satisfy (\ref{superliu0}) and then they are equal.

b) From (\ref{masterm}) and (\ref{masterbis})  it is now clear that for $z \in [z_c^- , z_c^+ ]$, the functions $j^{BST}(z, \cdot)$ and $j^{BIS}(z, \cdot)$ are  non-constant solutions of
\ben
\label{jfix}
{\bf J}(x) &=& \int_0 ^1 {\bf J}\left(xzu^{2z-1}\right) {\bf J}\left(xz(1-u)^{2z-1}\right) du\,,\cr
{\bf J}(0) &=& 1\,.
\een
With the same remark as above, we have uniqueness.

c) We see that (\ref{superliu}) yields (\ref{contraintedc}) and (\ref{contraintedc+}).
\QED
\medskip

To end this section let us notice that for the Yule-tree and the BST, 
the limit martingale connection (\ref{lmc1})  gives an important relation
 between the Laplace transforms:
\ben
\label{15}
j(z,x) = \int_0 ^\infty j^{BST}\left(z, x\frac{\eta^{2z-1}}{\Gamma(2z)}\right) \ e^{-\eta} d\eta\,.
\een

\section{Proof of theorems}
\subsection{Proof of Theorem \ref{bis=bst}}

The relation (\ref{a}) follows from (\ref{iterating}), (\ref{redefbis}), independence of $M_{\infty, u}^{BST}(z)$ with respect to $(U^{(v)} , v \prec u)$
 and the fact that $`EM_{\infty, u}^{BST}(z) = 1$ for every $u \in \BBu$.

To prove (\ref{b}) we first pass to the limit in $n$ to get a.s.
\begin{equation}
\label{above}
`E[M_\infty^{BST}(z) | {\cal F}_\infty^{BIS}] = M_\infty^{BIS}(z)\ .
\end{equation}
Set for a while, $X := M_\infty^{BST}(z)$, $Y := M_\infty^{BIS}(z)$ and ${\cal A} := {\cal F}_\infty^{BIS}$. 
Summarizing (\ref{above}) and  (\ref{lemme2}), we have
$$`E[ X | {\cal A}] = Y \ \ \hbox{and} \ \ \  X \buildrel{law}\over{=} Y\ .$$ 
From Exercise 1.11 in \cite{ChauYor} this implies $X=Y$ a.s.
\QED   

\subsection{Proof of Theorem \ref{MGEN=YULE}}
It is exactly the same line of argument as in the above subsection, using  (\ref{lemme1}) instead of (\ref{lemme2}) and (\ref{iter}) instead of (\ref{iterating}).

\subsection{Proof of Theorem \ref{derivMGEN=YULE}}\label{proof}

Since we work with fixed $z$, we omit it each time there is no possible confusion. 

The idea is to take advantage of the $L^1$ convergence of a multiplicative martingale and then come back taking logarithms. The proof can be given both in supercritical and in critical cases and we choose to present it for both cases, because it is not more complicated. Thus it will give an alternative proof of (\ref{limitmgen}) and the proof of (\ref{limitderivmgen}).

\medskip

\noindent{$\bullet$ First step: multiplicative martingales}. 
Notice that it could be more or less directly imported from theorem 3 in Kyprianou \cite{Kyp99} 
but ne- vertheless we give details
 to make the proof self-contained.

\medskip

Multiplicative martingales have appeared many times in the lite- rature, 
for instance in Neveu \cite{Nev1}, in Chauvin \cite{Chau91} in the branching Brownian motion framework, in Biggins and Kyprianou \cite{BK1} for discrete branching random walks and in Kyprianou \cite{Kyp99} for branching L\'evy processes. They are studied 
 in their own right
 in relation with functional equations or smoothing transformations and also, like here, 
 to help find information about additive martingales. Recall that here
\be
M^{GEN}_n(z)&=&\sum_{u \in {\cal L}_n} \  z^n e^{(1 -2z)S^{u}}\,,\\ 
M(t,z)&=& \sum_{u \in 
{\cal Z}_t} z^{|u|} e^{t(1 -2z)}\,.
\ee
Let for any real $y$,
\[
{\cal P}_n(y)= \prod_{|u| = n} j\big( y z^{n} e^{(1-2z) S^{u}} \big)\,,
\]
where $j(x) = j^{GEN}(z,x)$ (it is  a solution of equation (\ref{eqj}) 
with initial condition (\ref{ci})), and let
\[
{\cal P}(t)(y) = \prod_{u\in {\cal Z}_t} j\big( y z^{|u|} e^{(1-2z) t}\big)\,.
\]
To prove that ${\cal P}(t)(y)$ (resp. ${\cal P}_n(y)$) is a ${\cal F}_t$ (resp. ${\cal F}^n$)-martingale, decompose the set ${\cal Z}_t$ (resp. ${\cal L}_n$) with respect to a preceding line ${\cal Z}_s$ (resp. ${\cal L}_{n-1}$), apply the branching property and get the property of constant expectation: decompose $ `E({\cal P}(t)(y))$ (resp. $ `E({\cal P}_n(y))$) according to the first splitting time  $\tau_1$
and use the fact that $j$ is a solution of
equation (\ref{eqj}). 

Since $0\leq j\leq 1$, the martingale $\big({\cal P}(t)(y), t \geq 0\big)$ converges when $t$ goes to infinity 
a.s. and in $L^1$ to a limit ${\cal P}(\infty)(y)$; the martingale $\big({\cal P}_n(y), n \geq 0\big)$ converges when $n$ goes to infinity a.s. and in $L^1$ to a limit ${\cal P}_{\infty}(y)$.

\medskip

Let us now see why these two limits are equal: divide the set of individuals alive at time $t$ 
into those whose generation number is less than $n$ and those whose generation number is greater or equal to $n$. This gives a decomposition of ${\cal P}(t)(y)$; condition with respect to ${\cal F}^n$ and apply the branching property; fix $n$ and let $t$ tend to infinity to get 
\[
{\cal P}_n(y) = `E ({\cal P}(\infty)(y) | {\cal F}^n ) \ \ a.s.\,,
\]
and finally, since ${\cal P}(\infty)(y)$ is ${\cal F}^\infty :=\vee_n{\cal F}^n$ measurable
\ben
\label{newbis}
{\cal P}_{\infty}(y) ={\cal P}(\infty)(y) \ \ a.s.
\een

\noindent $\bullet$ Second step: back to the additive martingale, taking logarithms. 

\medskip
Let 
$(t_n)_{ n\geq 0}$ be a sequence going to infinity when $n$ goes to infinity.
We  use 
the behavior of 
 $j(z, .)$ near $0$.

Recall that for $z$ fixed, when $x\downarrow 0$, $\lim j(z,x) = 1$, so that  $-\log j(z,x)\sim 1-j(z,x)$. 
Moreover from (\ref{superliu0}) and (\ref{contraintedc}) we have the sharp estimates
\be
\lim_{x\rightarrow 0} \frac{1-j(z,x)}{x}
&=& 1 
\ \   \hbox{for} \ z \in (z_c ^- , z_c ^+)\,,
\\
\lim_{x\rightarrow 0} \frac{1-j(z,x)} {x |\log x| } &=& K_0 \ \ 
\hbox{for} \ z = z_c^\pm\,,
\ee
where $K_0 = 2/|2z-1|$. 
The quantity 
\[m_n(z) := \max \{ z^{|u|} e^{t_n (1 -2z)} \ ; \  u \in {\cal Z}_{t_n}\}\]
satisfies $m_n(z) \leq M( t_n,z)$. If $z = z_c^\pm$,  taking into account
that 
$\lim_n M( t_n,z_c^\pm ) = 0$, we have $\lim_n m_n (z_c ^\pm) = 0$ a.s. 
Now, for $z$ supercritical we check easily that
\begin{eqnarray*}
m_n(z) \leq m_n(z_c^+ )^{\log z/\log z_c^+}\ \ \ \hbox{if} \ 1 < z < z_c^+\,,\\
m_n(z) \leq m_n(z_c^-)^{\log z/\log z_c^-}\ \ \ \hbox{if} \ z_c^- \leq z < 1\,,
\end{eqnarray*}
and then $\lim_n m_n (z) = 0$ a.s.
Everything holds in the same way for the $n$-th generation and for $M^{GEN}_n(z)$ 
instead of $M(t_n,z)$.

We deduce that 
for every $\epsilon > 0$ there is some $n_0$ such that for every $n \geq n_0$ 
and $u \in {\cal Z}_{t_n}$ (resp. ${\cal L}_n$) 
\be
y (1- \epsilon) z^{|u|} e^{t_n (1 -2z)} \leq -\log j\big(z , y z^{|u|} e^{(1-2z) t_n} \big)\leq y (1+ \epsilon) z^{|u|}e^{t_n (1 -2z)}
\ee
in the supercritical case and
\be
(1- \epsilon)y (|u| \log z &+& t_n (1 -2z)) \, z^{|u|} e^{t_n (1 -2z)}\\ 
&+& K_0 (1- \epsilon) y (\log y) z^{|u|} e^{t_n (1 -2z)}\\
&\leq& -\log j\big(z , y z^{|u|} e^{(1-2z) t_n} \big) \\ 
&\leq& (1+ \epsilon)y (|u| \log z + t_n (1 -2z)) z^{|u|} e^{t_n (1 -2z)}\\
&+& K_0 (1+ \epsilon) y 
(\log y) z^{|u|} e^{t_n (1 -2z)}
\ee 
in the critical case
(resp. the same with $S^{u}$ instead of $t_n$). 
Adding up in $u$ we get respectively
\be
yM(t_n , z) (1- \epsilon) &\leq& -\log{\cal P}(t_n ,y)\leq y M(t_n , z) (1+ \epsilon)
\ee
and
\be
yM'(t_n , z) (1- \epsilon) &+& K_0 (y \log y) M(t_n , z)(1+ \epsilon)\\
&\leq& -\log{\cal P}(t_n ,y)\\
&\leq& yM'(t_n , z) (1+ \epsilon) + K_0 (y \log y) M(t_n , z)
\ee
(resp. the same -- mutatis mutandis -- for the GEN).
Taking limits in $n$, this implies a.s. 
\be -\log {\cal P}(\infty, y) = y M(\infty ,z)\,,\\
-\log {\cal P}_\infty (y) = y M_\infty^{GEN} (z)\,,
\ee
for $z$ supercritical and, for $z$ critical, the following
\be -\log {\cal P}(\infty, y) = y M'(\infty ,z)\,,\\
-\log {\cal P}_\infty (y) = y M_\infty^{'GEN} (z)\,.
\ee
With (\ref{newbis}), we now may conclude the proof. \QED

\subsection{Proof of Theorem \ref{cornew}} Let us consider only $z = z_c^+$ to  simplify.
We know that $j^{BST}(z, \cdot)$ and $j^{BIS}(z, \cdot)$ satisfy the same equation (\ref{jfix}).
Since we know by (\ref{contraintedc}) that $j^{BIS}(z, \cdot)$ satisfies
\ben
\label{liu}
\lim_{x \rightarrow 0+} \frac{1 - {\bf J}(x)}{x |\log x|} = \frac{2}{2z_c^+ - 1}\,,
\een
it is enough to prove that $j^{BST}(z, \cdot)$ satisfy also (\ref{liu}) (uniqueness mentioned in Section 4.2).
By Theorem \ref{derivMGEN=YULE} and (\ref{contraintedc}) again, $j(z, \cdot)$ satisfies also (\ref{liu}).
Now $j^{BST}(z, \cdot)$ is connected to $j(z, \cdot)$ by (\ref{15}). From  some elementary calculations and known properties of Laplace transforms we conclude that 
$j^{BST}(z, \cdot)$ 
satisfies also (\ref{liu}). 
\QED

\vspace{0.3cm}
\section{Links with Drmota's equations}\label{laplace}
 In this section, we make precise some probabilistic counterparts of solutions of equations introduced by M. Drmota in \cite{Dr1,Dr2,Dr0} as an analytical tool for a sharp study of the height of BST.

\subsection{General case}
We need some changes of parameter,  variables and  functions.
For $z \in [z_c^- , z_c^+]$ with $z\not= 1/2$, set 
\be
\alpha(z) = z^{\frac{1}{2z -1}}\,.
\ee
 For $z_c^- \leq z < 1/2$, (i.e. $c' \leq 2z < 1$) the function $\alpha$ 
increases from $\alpha_{c'}=e^{1/c'}$ to $+\infty$, and when $1/2 < z \leq z_c^+$ (i.e. $1 < 2z \leq c$) it increases from $0$ to 
$\alpha_c = e^{1/c}$. We  often write  $\alpha$ instead of $\alpha(z)$ to simplify.

For $2z \not= 1$ let
$\varphi(z, x) = x^{-1}j(z, x^{1-2z})$. 
Equation (\ref{eqj}) is translated into
\ben
\label{11}
\varphi(z, x) &=& \alpha^{-2} \int_x ^\infty \varphi(z , y/\alpha)^2 dy\,.
\een
Moreover, equation (\ref{ci}) becomes:
\ben
\label{cim+}
\lim_{x \rightarrow \infty} x \varphi (z, x) &=& 1 \ \  , \ \hbox{if}\ \ \ 2z > 1\,,\\
\label{cim-}
\lim_{x \rightarrow 0} x \varphi (z, x) &=& 1 \ \  , \ \hbox{if}\ \ \ 2z < 1\,.
\een
and equation $(\ref{superliu0})$ becomes:
\ben
\label{contrainted+}
\lim_{x\rightarrow\infty}\frac{1 - x\varphi(z,x)}{x^{1-2z}} &=&  1 \ \  , 
\ \hbox{if}\ \ \  1 < 2z < c\,,\\
\label{contrainted-}
\lim_{x\rightarrow0}\frac{1 - x\varphi(z,x)}{x^{1-2z}} &=& 1  \ \  , \ \hbox{if}\ \ \ c' < 2z < 1\,.
\een

Drmota used the solution of the retarded differential equation
\ben
\label{phi}
\Phi'(x)&=& -\frac{1}{\alpha^2} \ \Phi(x/\alpha)^2\,,\\
\label{ciphi}
\Phi(0) &=& 1\,,
\een
and also the solution of the (retarded) convolution equation
\ben
\label{psi}
y\Psi(y/\alpha) = \int_0 ^y \Psi(w)\Psi(y -w) dw\,.
\een
If $\varphi$ satisfies (\ref{phi}), then $\varphi^\kappa$ 
defined by $\varphi^\kappa (x) = \kappa \varphi(\kappa x)$ satisfies the same equation.
Similarly, if $\psi$ satisfies (\ref{psi}), then $\psi^\kappa$ defined by $\psi^\kappa (u) =  \psi(u/\kappa )$ satisfies the same equation.
Drmota (Lemmas 18, 19 and 23 of \cite{Dr1} and Prop. 5.1 of \cite{Dr0})  proved that for $ 1 < \alpha \leq \alpha_c$ 
, there is  a unique entire solution $\varphi_\alpha$ of (\ref{phi}), and that $\varphi_\alpha$ is a Laplace transform: 
\ben
\label {10}
\varphi_\alpha (x) = \int_0 ^{\infty} \psi_\alpha (y) e^{-x y} dy\,,
\een
where $\psi_\alpha$ is solution of (\ref{psi}). Moreover, these functions have the following behavior:
\ben
\label{varphi}
\lim_{x \rightarrow \infty}\frac{1 - x\varphi_\alpha(x)}{x^{1 -2z}} &=& K_1\,, \ \ \hbox{if}\ \ \ 1 < \alpha < \alpha_c\,,\\
\label{varphic}
\lim_{x \rightarrow \infty}\frac{1 - x\varphi_\alpha(x)}{x^{1 -2z}\log x} &=& K_2\,, \ \ \hbox{if}\ \ \ \alpha = \alpha_c\,,
\een
where   $K_1$ and $K_2 \in (0, \infty)$.
This implies in particular that $\varphi_\alpha$ satisfies (\ref{11}). But, if $\varphi$ satisfies (\ref{phi}), then $\varphi^\kappa$ defined by $\varphi^\kappa (x) = \kappa \varphi(\kappa x)$ satisfies the same equation, so that the role of the initial condition (\ref{ciphi}) is to fix the constant $\kappa$. Conversely, we may 
choose $\kappa(z)$
 such that $\varphi_{\alpha(z)}^{\kappa(z)}$ satisfies  (\ref{varphi}) (resp. (\ref{varphic})) 
 with  limit $1$ instead of $K_1$ (resp. $K_2$). Comparing with (\ref{contrainted+}), we have (at least for $1 < z < z_c^+$) 
 that $\varphi_{\alpha (z)}^{\kappa(z)} = \varphi (z, \cdot)$. 
 
For $1 < z < z_c^+$, we know  from  (\ref{lemme1}) and Theorem \ref{MGEN=YULE} that  the Yule-time limit martingale $M(\infty , z)$ and the Yule-generation limit martingale $M^{GEN}_\infty(z)$ 
 have the same Laplace transform. We see now that  
\ben
\label{dr1}
j(z,x)= j^{GEN}(z,x) = x^{\frac{1}{1-2z}} \kappa\varphi_\alpha 
\big( \kappa x^{\frac{1}{1-2z}}\big)\,,
\een
or equivalently
\ben
\label{dr2}
`E e^{-x^{1-2z}M(\infty ,z)} = `E e^{-x^{1-2z} M^{GEN}_\infty (z)} = \kappa x \varphi_\alpha (\kappa x)\,.
\een
for some constant $\kappa > 0$.

In the critical case, 
 $z= z_c^+$,  we use the derivative martingales. From Theorem \ref{derivMGEN=YULE} and Section 4.2  
 we see that (\ref{dr1}) holds again or, in other words,
\ben
\label{dr3}
`E e^{-x^{1-2z_c^+}|M'(\infty ,z_c^+)|} = `E e^{-x^{1-2z_c^+} |M^{'GEN}_\infty (z_c^+)|} = \kappa x \varphi_{\alpha_c} (\kappa x)\,.
\een
for some $\kappa$.

Now let us consider
 $\psi_\alpha$ 
as defined in 
(\ref{psi})
 or (\ref{10}). The relation (\ref{15}) together with (\ref{10}) and (\ref{dr1}) gives easily
$$\psi_\alpha(y/\kappa) = j^{BST}\left(z, \frac{y^{2z-1}}{\Gamma(2z)}\right) \,.$$
 The BST limit martingale $M^{BST}_\infty (z)$ and the Bisection limit martingale $M^{BIS}_\infty(z)$ 
(which are equal in the supercritical case) have the same Laplace transform 
 and 
$$`Ee^{-\frac{y^{2z -1}}{\Gamma (2z)}M^{BST}_\infty (z)} = `Ee^{-\frac{y^{2z -1}}{\Gamma (2z)}M^{BIS}_\infty (z)}= \psi_\alpha (y/\kappa)\,.$$ 
In the critical case, $z=z_c^+$, 
 we have (recall Theorem \ref{cornew})
$$`Ee^{-\frac{y^{2z_c^+ -1}}{\Gamma (2z_c^+)}|M'^{BST}_\infty (z_c^+)|} = `Ee^{-\frac{y^{2z_c^+ -1}}{\Gamma (2z_c^+)}|M'^{BIS}_\infty (z_c^+)|}=\psi_{\alpha_c} (y/\kappa)\,.$$ 

\vspace{0.3cm}
\subsection{Particular cases}

1) For $z = 1/2$, we have $M(t, 1/2) \equiv 1$ and $j(1/2 , x) = e^{-x}$. In this case $\alpha$ is not defined.

\smallskip

2) For $z = 1$, as previously mentioned, we have $M(t, 1) = e^{-t} N_t$ 
whose limit is $\xi \sim {\cal E}(1)$, of Laplace transform $j(1, x) = \frac{1}{1+x}$. 
This yields $\varphi(1, x) = \frac{1}{1+x}$, which corresponds of course 
to $\varphi_\alpha (x)$ with $\alpha = \alpha (1) = 1$.

\smallskip

3) The function $\Phi(x) = e^{-x/4}$ is solution of (\ref{phi}) and (\ref{ciphi}) for $\alpha=2$, but it does not correspond to the limit of a martingale from our families. In fact, there is no real $z \in (0, \infty)$ such that $\alpha (z) = 2$.

\smallskip

4) 
Drmota \cite{Dr1} in his Lemma 11 noticed that, for $\alpha = 16$ 
$$\bar\varphi (x) = \frac{1 + x^{1/4}}{x}\, e^{-x^{1/4}}$$
satisfies equation (\ref{phi}) and the initial condition (\ref{cim-}).

In $(0, \infty)$ the equation $\alpha(z) = 16$ has two solutions: $z = 1/8$ and $z = 1/4$. Since $1/8 < z_c ^-$, the corresponding martingale $M(t, 1/8)$ goes to $0$ a.s. and then $\varphi (1/8, x) \equiv 1$. Besides, we have $z_c^- < 1/4 < 1/2$, so that $j(1/4 , x)$ is not constant and satisfies (\ref{eqj}), (\ref{ci}), (\ref{superliu0}). Again, 
 for every constant $\kappa$, the function $\bar\varphi^\kappa =\kappa\bar\varphi(\kappa .)$  satisfies the same system. The constraint 
(\ref{contrainted-}) leads to $\kappa = 4$.
Taking into account the uniqueness, 
we see that
\be`E \exp -xM(\infty , 1/4) = 
(1 + \sqrt{2 x}) e^{-\sqrt{2 x}}\,,
\ee 
which allows identification of the law of $M(\infty , 1/4)$.
As it is noticed in another context in  Yor \cite{yorexpo} pp. 110--111, this is the Laplace transform of the density
\be
G(t) = \frac{1}{\sqrt{2\pi t^5}} \ e^{-1/2t} \ \ , \ \ t\geq 0\,,
\ee
of $\big(2\gamma_{3/2})^{-1}$, where we use the notation $\gamma_\alpha$ for a variable with distribution gamma with parameter $\alpha$, i.e.
$$`P(\gamma_\alpha \in dy) = \frac{1}{\Gamma(\alpha)}y^{\alpha -1} e^{-y} dy\ \ \ \ \  y > 0 \,.$$
In other words
\ben
\label{mdef}
M(\infty, 1/4) \buildrel{law}\over{=} \big(2\gamma_{3/2})^{-1}\,.
\een
Let us now look for the law of $M_\infty^{BST}(1/4)$.
It is easier to work with $A := \big(M_\infty^{BST}(1/4)\big)^{-2}$. 
The limit connection (\ref{lmc1}) yields
\be
\big(2\gamma_{3/2}\big)^2  \buildrel{law}\over{=} \pi \gamma_1 \cdot A\,,
\ee
where $A$ is independent of $\gamma_1$ (which is our $\xi$). 
Introducing the Mellin transform we have
\ben
`E\big(2\gamma_{3/2}\big)^{2s} = \pi^s `E (\gamma_1)^s `E (A^s)\,.
\een
Recall that $`E (\gamma_a)^s = \frac{\Gamma(s+ a)}{\Gamma(a)}$ for $\Re s > -a$, 
so that, for $\Re s > -3/4$,
\ben
\label{emus}
`E (A^s) = \frac{2^{2s} \Gamma(2s+ \frac{3}{2})}{\Gamma(\frac{3}{2})} 
\ \frac{1}{\pi^s \Gamma (s+ 1)}\,.
\een
With the duplication formula 
\ben
\label{dupli1}
\Gamma(2y) = (2\pi)^{-1/2} 2^{2y -\frac{1}{2}} \Gamma(y) \Gamma(y + \frac{1}{2})
\een
taken at $y = s + \frac{3}{4}$ and $y = \frac{3}{4}$,
the identity (\ref{emus}) becomes
\ben
\label{newmus}
`E (A^s) = \left(\frac{16}{\pi}\right)^s 
\frac{\Gamma (s + \frac{3}{4})}{\Gamma(\frac{3}{4})\Gamma (s+1)} 
\frac{\Gamma (s + \frac{5}{4})}{\Gamma ( \frac{5}{4})}\,.
\een
Now, for $a, b > 0$, a beta random variable $\beta_{a,b}$
with density
$$`P(\beta_{a,b} \in dy) = \frac{\Gamma(a+b)}{\Gamma(a) \Gamma(b)} y^{a-1} (1-y)^{b-1} 
\ \ \ \,,  \ \ y \in (0,1)\,,$$
has Mellin transform:
$$`E (\beta_{a,b})^s = \frac{\Gamma(a+s)\Gamma(a+b)}{\Gamma(a)\Gamma (a+b+s)}\,,\ \ \ \ \Re s > -a\,,$$ which allows to transform (\ref{newmus}) into
\be
`E (A^s) = \left(\frac{16}{\pi}\right)^s `E(\beta_{3/4, 1/4})^s \ `E(\gamma_{5/4})^s\,,\ \ \ \ \Re s > -3/4\,,
\ee
and then $
A\buildrel{law}\over{=} \frac{16}{\pi} \ \beta_{3/4, 1/4} \cdot \gamma_{5/4}
$ or , coming back to $M^{BST}$
\ben
\label{bstdef}
M^{BST}_\infty (1/4) \buildrel{law}\over{=} \frac{\sqrt{\pi}}{4}\big(\beta_{3/4, 1/4}\big)^{-1/2} \cdot \big(\gamma_{5/4}\big)^{-1/2}\,,
\een
where the two variables in the right-hand side are independent.

\begin{rem}
\end{rem}
The limit martingale connection (3) rewritten with (\ref{mdef}) and (\ref{bstdef}) gives
\ben
\label{yor}
\big(\gamma_{3/2}\big)^2 \buildrel{law}\over{=} 4 \beta_{3/4, 1/4} . \gamma_1 . \gamma_{5/4}\,.
\een
This identity can be viewed as a consequence of the two following identities in law that can be found in Chaumont-Yor \cite{ChauYor}. The first one is, for any parameter $a$
\ben
\big(\gamma_a\big)^2 \buildrel{law}\over{=} 4 \gamma_{\frac{a}{2}} . \gamma_{\frac{a+ 1}{2}}\,,
\een 
which is the stochastic version of (\ref{dupli1}), and the second one (coming from the classical beta-gamma algebra) is
\ben
\gamma_{\frac{a}{2}} \buildrel{law}\over{=} \beta_{\frac{a}{2}, 1 - \frac{a}{2}} . \gamma_1\,,
\een
holding for $a < 2$. Both give
\ben
\big(\gamma_a\big)^2  \buildrel{law}\over{=} 4 \beta_{\frac{a}{2}, 1 - \frac{a}{2}} . \gamma_1 . \gamma_{\frac{a+ 1}{2}}\,.
\een 
Identity (\ref{yor}) is obtained taking $a = 3/2$.

\bigskip

\noindent{\bf Acknowledgements}
\bigskip

We thank Michael Drmota for many valuable discussions on the respective merits 
of the martingale approach and the analytic approach, and Marc Yor for pointing out the above remark.

\bibliographystyle{plain} 
\small
\bibliography{chaurouarxiv1}
\end{document}